\documentclass{amsart}

\usepackage{amssymb}
\usepackage[all]{xy}
\usepackage{graphics}

\newtheorem*{introthm}{Theorem}

\newtheorem{theorem}{Theorem}[section]
\newtheorem{lemma}[theorem]{Lemma}
\newtheorem{proposition}[theorem]{Proposition}
\newtheorem{corollary}[theorem]{Corollary}

\theoremstyle{definition}
\newtheorem{definition}[theorem]{Definition}
\newtheorem{example}[theorem]{Example}
\newtheorem{remark}[theorem]{Remark}
\newtheorem{problem}[theorem]{Problem}

\def\Chi{{\rm X}}

\def\div{{\rm div}}

\def\tdr{{\rm tdr}}
\def\quot{/\!\!/}
\def\mal{\! \cdot \!}

\def\tq{\mathbin{\!{/ \! \! \lower 3pt \hbox{$\scriptscriptstyle{\rm
          tq}$}}\!}}

\def\rq#1{\widehat{#1}}
\def\t#1{\widetilde{#1}}
\def\b#1{\overline{#1}}

\def\CC{{\mathbb C}}
\def\KK{{\mathbb K}}

\def\ZZ{{\mathbb Z}}

\def\QQ{{\mathbb Q}}
\def\PP{{\mathbb P}}

\def\h#1{\widehat{#1}}

\def\id{{\rm id}}
\def\CDiv{{\rm CDiv}}

\def\codim{{\rm codim}}
\def\Supp{{\rm Supp}}
\def\Spec{{\rm Spec}}
\def\conv{{\rm conv}}
\def\cone{{\rm cone}}

\def\pr{{\rm pr}}

\begin{document}



\title[Quotients of divisorial toric varieties]
      {Quotients of divisorial toric varieties}

\author[A.~A'Campo-Neuen]{Annette A'Campo-Neuen} 
\address{Fachbereich Mathematik, Johannes Gutenberg-Universit\"at Mainz,
  55099 Mainz, Germany}
\email{acampo@mathematik.uni-mainz.de}

\author[J.~Hausen]{J\"urgen Hausen} 
\address{Fachbereich Mathematik und Statistik, Universit\"at Konstanz,
  78457 Konstanz, Germany}
\email{Juergen.Hausen@uni-konstanz.de}


\begin{abstract}
We consider subtorus actions on divisorial toric varieties.
Here divisoriality means that the variety has many Cartier divisors
like quasiprojective and smooth ones. We characterize when a subtorus 
action on such a toric variety admits a categorical quotient in
the category of divisorial varieties. Our result generalizes previous
statements for the quasiprojective case.
An important tool for the proof is a universal reduction of an
arbitrary toric variety to a divisorial one. This is done in terms of
support maps, a notion generalizing support functions on a polytopal
fan. A further essential step is the decomposition of a given subtorus
invariant regular map to a divisorial variety  into an invariant toric
part followed by a non-toric part.
\end{abstract}

\maketitle

\section*{Introduction}

It is a frequently occuring question in algebraic geometry, if an
algebraic group action $G \times X \to X$ admits a categorical
quotient, i.e., a regular map $X \to Y$ that is universal with respect
to $G$-invariant regular maps $X \to Z$. For example, moduli
functors are often corepresented by categorical quotients. 
In general, it is a difficult problem to decide whether a
categorical quotient exists. Some counterexamples for actions of the
multiplicative group $\CC^{*}$ are presented in \cite{acha4}.

As these examples show, difficulties already arise with
subtorus actions on toric varieties. Such actions have been investigated
by several authors, mainly focusing on the much more restrictive
concept of a good quotient, see e.g. \cite{KaStZe}, \cite{Sw} and
\cite{Hm}. The description of toric varieties in terms of rational
fans relates the problem of constructing quotients to problems of
combinatorial convexity. Hence the class of toric varieties serves as
a testing ground for more general ideas.

Now, let $X$ be a toric variety and let $H$ be a subtorus of the big
torus of $X$. Our approach to categorical quotients for the induced action
of $H$ on $X$ is to consider the problem in suitable subcategories. A
first step is to construct a quotient in the category of toric
varieties itself: In \cite{acha1}, we showed that there always exists
a {\it toric quotient\/}
$$p \colon X \to X \tq H .$$

This is a toric morphism that is universal with respect to
$H$-invariant toric morphisms. The essential part of the proof is
an explicit algorithm in terms of combinatorial data. The toric
quotient is a canonical starting point for quotients in further
categories. For example, in \cite{acha2} we gave an explicit method
to decide by means of the toric quotient when a subtorus
action on a quasiprojective toric variety admits a
categorical quotient in the category of quasiprojective varieties.

In the present article we give a considerable generalization of
the results of \cite{acha2}, namely we solve the analogous problem
in the category of divisorial varieties. Recall that an irreducible
variety $X$ is called divisorial if every point $x \in X$ has an
affine neighbourhood of the form $X \setminus \Supp(D)$ with an effective
Cartier divisor $D$ on $X$, see e.g. \cite{Bo} and \cite[II.2.2]{SGA}.

The class of divisorial varieties contains the quasiprojective
varieties as well as all $\QQ$-factorial varieties. It has nice
functorial properties, see \cite{Bo}, and moreover it often provides a
natural framework to extend statements known to hold for quasiprojective
varieties on the one hand and for smooth varieties on the other hand.

A connection to toric geometry is provided by the embedding results of
\cite{ha}: A variety is divisorial if and only if it admits a closed
embedding into a smooth toric prevariety $Z$ having an affine diagonal
map $Z \to Z \times Z$. The equivariant version of this statement
implies in particular that a toric variety is divisorial if and only
if it has enough invariant effective Cartier divisors in the sense of
T.~Kajiwara \cite{Ka}, see Section~1.

Now, given a divisorial toric variety $X$ and a subtorus $H$ of the
big torus of $X$, when does the action of $H$ on $X$ admit a
categorical quotient in the category of divisorial varieties? 
As mentioned, we start with the toric quotient
$$p \colon X \to X \tq H.$$

A first problem is that in general the toric quotient space $X \tq H$
is not a divisorial variety. To deal with this effect, we construct a {\it
toric divisorial reduction\/}. This is a toric morphism 
$$q \colon X \tq H \to (X \tq H)^{\tdr}$$
which is universal with respect to toric morphism to divisorial toric
varieties. The question then is, how these toric constructions behave
in the essentially larger category of arbitrary divisorial varieties. Our
main result gives the following answer, see
Corollary~\ref{divcatquot}:

\begin{introthm}
The action of $H$ on $X$ admits a categorical quotient in the category
of divisorial varieties if and only if the composition $q \circ
p$ is surjective. Moreover, in the latter case, $q \circ p$ is the
desired categorical quotient.
\end{introthm}

The paper is organized as follows: In Section~\ref{section1} we
discuss divisoriality in the context of $G$-varieties and provide some
general statements used in the subsequent constructions. 
Sections~\ref{section2} and~\ref{section3} are devoted to the
construction of the toric divisorial reduction. This is done in the
language of combinatorial convexity. The main tool are convex support
maps extending the notion of a convex support function on a fan.

Generalizing the corresponding well-known statement on projectivity
and support functions, we show that divisoriality of a given toric
variety is characterized by the existence of a strictly convex support
map on its fan. Moreover, we relate convex support maps to toric
morphisms to divisorial toric varieties. This allows the
construction of the toric divisorial reduction. Finally, we
present some examples in Section~\ref{section3}.

In Sections~\ref{section4} and~\ref{section5} we prepare the proof of
the main results. The essential task is to reduce arbitrary
$H$-invariant regular maps to $H$-invariant toric morphisms. This is
done by the Decomposition Lemma presented in Section~\ref{section5}:
Given an $H$-invariant regular map $f \colon X \to Y$ to a divisorial
variety, we construct a decomposition $f = h \circ g$ with an
$H$-invariant toric morphism $g$ followed by a rational map $h$
defined near $g(X)$.  

The ingredients for the proof of this Decomposition Lemma are the
abovementioned embedding of $Y$ into a certain smooth toric prevariety
$Z$ provided by~\cite{ha} and the following lifting result, presented
in Section~\ref{section4}: There exist quasiaffine toric varieties
$\t{X}$ and $\t{Z}$ ``above'' $X$ and $Z$ respectively such that the
map $f$ admits a lifting $\t{f} \colon \t{X} \to \t{Z}$. This basically
reduces the decomposition problem to the case of quasiaffine toric
varieties.

In Section~\ref{section6} we give statements and proofs of the main
results. Finally, in Section~\ref{section7} we formulate an open
problem on categorical quotients for subtorus actions on toric
varieties.

\section{Divisorial $G$-varieties}\label{section1}

Throughout the whole article, we work over a fixed
algebraically closed field $\KK$. So a prevariety is a reduced
irreducible scheme of finite type over $\KK$, and a variety is a 
separated prevariety. We say that a prevariety $X$ is of 
{\it affine intersection\/}, if its diagonal morphism $X \to X \times X$
is affine. 

As usual, when we speak of a $G$-(pre-)variety where $G$ is
an algebraic group, we mean an algebraic (pre-)variety $X$ together
with a $G$-action given by a regular map $G \times X \to X$.
For the basic notions on toric varieties and prevarieties, we
refer to~\cite{Fu} and~\cite{acha3}.

In this section, we provide some general facts on group actions on
divisorial varieties. Following Borelli \cite{Bo}, we call a 
prevariety $X$ {\it divisorial\/} if every point $x \in X$ has an
 affine open neighbourhood of the form $U = X \setminus \Supp(D)$ with
an effective Cartier divisor $D$ on $X$. 

\begin{remark}\label{divisorialproperties}
\begin{enumerate}
\item Quasiprojective varieties are divisorial.
\item Locally closed subspaces of divisorial prevarieties are
  divisorial. 
\item Every divisorial prevariety $X$ is of affine intersection.
\item Every $\QQ$-factorial prevariety of affine intersection is
  divisorial.
\end{enumerate}
\end{remark}

A {\it geometric quotient\/} for the action of a reductive group $G$ on
a variety $X$ is an affine regular map $p \colon X \to Y$
such that the fibres of $p$ are precisely the $G$-orbits and the
canonical homomorphism $\mathcal{O}_{Y} \to p_{*}(\mathcal{O}_{X})^{G}$
is bijective. The analogous notion in the setting of prevarieties,
i.e. for possibly non-separated $X$ and $Y$, is called a 
{\it geometric prequotient\/}. 

In the sequel, we shall make use of the following characterization of
divisoriality in terms of geometric quotients and closed embeddings,
see~\cite[Theorem~3.1]{ha}: 

\begin{theorem}\label{einbettungssatz}
A variety $X$ is divisorial if and only if one of the following
 statements holds:
\begin{enumerate}
\item $X$ is a geometric quotient of a quasiaffine variety by a free
  algebraic torus action.
\item $X$ admits a closed embedding into a smooth toric prevariety of
  affine intersection.
\end{enumerate}
\end{theorem}

Here a torus action is called free if every orbit map is a locally
closed embedding. The above result has the following equivariant version,
see~\cite[Theorem~3.4]{ha}:

\begin{theorem}\label{aequiveinbettung}
Let $X$ be a normal divisorial $T$-variety where $T$ is an algebraic
torus acting effectively.
\begin{enumerate}
\item There is a quasiaffine variety $\rq{X}$ with a regular action of
  a torus $T \times H$ such that $H$ acts freely with a
  $T$-equivariant geometric quotient $\rq{X} \to X$.
\item There is a $T$-equivariant closed embedding $X \to Z$ into a
  smooth toric prevariety $Z$ of affine intersection where $T$ acts as
  a subtorus of the big torus.   
\end{enumerate}
\end{theorem}

A first consequence is that divisorial varieties with torus actions
always have many invariant effective Cartier divisors. For
a toric variety this means that it is divisorial if and only if
it has {\it enough invariant effective Cartier divisors\/}  in the sense
defined by T.~Kajiwara, see \cite{Ka}.

\begin{proposition}\label{enoughinvar}
Let $T$ be an algebraic torus, and let $X$ be a normal algebraic
$T$-variety $X$. Then $X$ is
divisorial if and only if there exist $T$-invariant effective Cartier
divisors $D_{1}, \ldots, D_{r}$ on $X$ such that the sets $X \setminus
\Supp(D_{i})$ are affine and cover $X$.
\end{proposition}

\proof We may assume that $T$ acts effectively. Let $X$ be
divisorial. By Theorem~\ref{aequiveinbettung}, there is a 
$T$-equivariant closed embedding of $X$ into a smooth toric prevariety
$Z$ of affine intersection where $T$ acts as a subtorus of the big
torus. Hence $X$ inherits the desired property from $Z$. The reverse
implication is trivial. 
\endproof

As the example of the rational nodal curve with standard
$\KK^{*}$-action shows, the assumption of normality is essential in
the above statement. Our next result states that divisoriality is
inherited by geometric quotients for torus actions:

\begin{proposition}\label{div2quot}
Let $T$ be an algebraic torus and suppose that $X$ is a normal
$T$-variety with geometric quotient $p \colon X \to Y$. Then $X$
is divisorial if and only if $Y$ is divisorial. 
\end{proposition}

\proof We may assume that the torus $T$ acts effectively on $X$. If
the quotient variety $Y$ is divisorial, then we obtain the desired
effective Cartier divisors on $X$ by pulling back suitable divisors
from $Y$. Conversely, suppose that $X$ is divisorial. Then,
by Theorem~\ref{aequiveinbettung}, we may
assume in the proof that $X$ is a quasiaffine $T$-variety. 

Given $y \in Y$, we have to find an affine open neighbourhood of
$y$ that is the complement of the support of an effective Cartier
divisor on $Y$. Choosing any $T$-equivariant affine closure of $X$, we
find a function $f \in \mathcal{O}(X)$, homogeneous with respect to
some character $\chi_{f} \in \Chi(T)$, such that for $D := \div(f)$
the $T$-invariant set $U := X\setminus V(f)=X \setminus \Supp(D)$ is
an affine neighbourhood of the fibre $p^{-1}(y)$. 

By $T$-closedness of $p \colon X \to Y$, the set $V := p(U)$ is an
open neighbourhood of $y \in Y$. Moreover, as a geometric quotient
space of the affine $T$-variety $U$, the set $V$ is again affine.
Thus, to prove the assertion, we only have to show that $p(\Supp(D))$
is the support of an effective Cartier divisor $E$ on $Y$. We
construct local equations for such an $E$. 

First we claim that every point $z \in Y$ has an affine neighbourhood $V_z
\subset Y$ such that on $U_z := p^{-1}(V_z)$ there is an invertible
function $h_z \in \mathcal{O}(U_z)$ that is homogeneous with respect
to some positive multiple $m_z \chi_{f}$. To check this, start with any affine
neighbourhood $V_z \subset Y$ of $z$ and choose a point $x \in p^{-1}(z)$. 
Consider the sublattice $N \subset \Chi(T)$ of  characters occuring as
weights of homogeneous functions $g \in \mathcal{O}(U_z)$ with $g(x) =  1$.

The sublattice $N$ is of full rank in $\Chi(T)$: Otherwise we found a
nontrivial one-parameter-subgroup $\lambda \colon \KK^{*} \to T$ such
that $\chi \circ \lambda = 1$ holds for all $\chi \in N$. It 
follows that $\lambda(\KK^{*})$ is contained in the isotropy group
$T_{x}$. On the other hand, the $T$-action on $U_z$ is effective and
closed. Hence $T_{x}$ is finite, a contradiction. 
Thus $N$ is of full rank. In particular, some positive
multiple $m_z \chi_{f}$ lies in $N$ and our claim follows.

Now cover $Y$ by finitely many $V_z$ as in the above claim. Then we
may assume that all the invertible functions $h_z \in
\mathcal{O}(U_z)$ are homogeneous with respect to the same multiple $m
\chi_f$. Every function $g_z := f^m/h_z$ is $T$-invariant, regular on
$U_z$ and vanishes precisely on $\Supp(D) \cap U_z$. Since it is
$T$-invariant, $g_z$ may be viewed as a regular function
on $V_z = p(U_z)$ and there its zero set is just
$$p(\Supp(D) \cap U_z) = p(\Supp(D)) \cap V_z. $$
Since every $g_z / g_{z'}$ is an invertible regular function on $V_z
\cap V_{z'}$ it follows that the $g_z$ are local equations for the
desired Cartier divisor $E$ on $Y$. \endproof

As T.~Kajiwara has shown, every toric variety $X$ with enough
invariant effective Cartier divisors arises as a geometric quotient of
a quasiaffine toric variety $\rq{X}$ by an algebraic subgroup of the
big torus of $\rq{X}$, see~\cite[Theorem~1.9]{Ka}. In view of the
above results, we can enhance Kajiwara's statement as follows:

\begin{corollary}\label{divtorvar}
A toric variety $X$ is divisorial if and only if there is a
quasiaffine toric variety $\rq{X}$ and a toric morphism $p \colon
\rq{X} \to X$ such that $\ker(p)$ is a subtorus of the big torus of
$\rq{X}$ and $p$ is a geometric quotient for the action of $\ker(p)$
on $\rq{X}$.
\end{corollary}

\proof If $X$ is divisorial, then Theorem~\ref{aequiveinbettung} gives
the desired quotient presentation. The converse follows from 
Proposition~\ref{div2quot}. \endproof 

Finally, we consider translates of divisorial open subsets with respect
to an action of a connected group. If the complement of the subset
is small enough, the union of such translates is again divisorial: 

\begin{lemma}\label{translates}
Let $G$ be a connected linear algebraic group, and let $X$ be a normal
$G$-variety. If $U \subset X$ is a divisorial open subset with
$\codim(X \setminus U) \ge 2$, then also $G \mal U$ is divisorial.
\end{lemma}

\proof We may assume that $X = G \mal U$ holds. Let $D_{1}^{U}, \ldots,
D_{r}^{U}$ be Cartier divisors on $U$ such that the sets $U_{i} := U
\setminus \Supp(D_{i}^{U})$ form an affine cover of $U$. By closing
components, each $D_{i}^{U}$ extends to a Weil divisor $D_{i}$ on $X$.

We claim that $X \setminus \Supp_{D_{i}} = U_{i}$. To see this, let
$A_{i} := X \setminus U_{i}$. Since $U_{i}$ is affine, $A_{i}$ is of
pure codimension one. Clearly $\Supp(D^{U}_{i}) \subset A_{i}$ and
hence $\Supp(D_{i}) \subset A_{i}$. Thus $\Supp(D_{i})$ is a union of
irreducible components of $A_{i}$. Moreover we have
$$ X \setminus U = X \setminus (U_{i} \cup \Supp(D^{U}_{i})) = A_{i}
\setminus \Supp(D^{U}_{i}). $$

Since $X \setminus U$ has codimension at least two, it follows that
for each irreducible component $A_{i}'$ of $A_{i}$ its intersection
with $\Supp(D^{U}_{i})$ is dense in $A_{i}'$. This implies $A_{i} =
\Supp(D_{i})$ and our claim is proved. In particular, we have
$$ X = G \mal U = G \mal \bigcup_{i = 1}^{r} X \setminus
\Supp(D_{i}) = \bigcup_{i=1}^{r} \bigcup_{g \in G} X \setminus \Supp(g
\mal D_{i}). $$

Thus it suffices to show that for each $D_{i}$ some multiple is
Cartier on $X$. This is done as follows: The restriction $D'_{i}$ of
$D_{i}$ to the regular locus $X' \subset X$ is Cartier. Since $X'$ is
$G$-invariant, we may apply $G$-linearization, i.e., replacing $D_{i}$
with a suitable multiple we achieve that $\mathcal{O}_{D'_{i}}$ is a
$G$-sheaf, see e.g.~\cite[Proposition~2.4]{DMV}.

We claim that this structure of a $G$-sheaf extends canonically to
$\mathcal{O}_{D_{i}}$. For an open set $V
\subset X$ let $V' := V \cap X'$. Given a section $s \in
\mathcal{O}_{D_{i}}(V)$, we define its translates $g \mal s$ as follows:
Translate the restriction $s' \in \mathcal{O}_{D_{i}}(V')$ to a section
$g \mal s' \in \mathcal{O}_{D_{i}}(g \mal V')$ and then extend $g \mal s'$
to the desired section $g \mal s \in \mathcal{O}_{D_{i}}(g \mal V)$.

Using the $G$-sheaf structure on $\mathcal{O}_{D_{i}}$ we see that locally
$\mathcal{O}_{D_{i}}$ is generated by a single function. That means $D_{i}$
 is a Cartier divisor. \endproof

\section{Support maps}\label{section2}

Projectivity of a given toric variety is characterized by
the existence of a strictly convex support function on its fan, see
e.g.~\cite{Fu}. Generalizing the notion of a support function
here we introduce the concept of a support map on a fan and 
define convexity properties for such maps. The main
result of this section states that for a given fan existence of a
strictly convex support map is equivalent to divisoriality of the
associated toric variety. 

For a lattice $N$, we denote the associated rational vector
space by $N_{\QQ}$. A {\it cone\/} in $N$ is a polyhedral (not
necessarily strictly) convex cone $\sigma \subset N_{\QQ}$.
A {\it quasifan\/} in $N$ is a finite set $\Lambda$ of cones
in $N$ such that for $\sigma \in \Lambda$ also every face of $\sigma$
belongs to $\Lambda$ and for $\sigma, \sigma' \in \Lambda$ the
intersection $\sigma \cap \sigma'$ is a face of both, $\sigma$ and
$\sigma'$. A {\it fan\/} is a quasifan containing only strictly convex
cones.

The {\it support\/} of a quasifan $\Lambda$ is the union of all
its cones and is denoted by $\vert \Lambda \vert$. A {\it map of
  quasifans\/} $\Lambda$ in a lattice $N$ and $\Lambda'$ in a lattice
$N'$ is a lattice homomorphism $N \to N'$ such that the associated
linear map $N_{\QQ} \to N'_{\QQ}$ maps the cones of $\Lambda$ into
cones of $\Lambda'$.

For the definition of support maps, fix a lattice $N$ and a quasifan
$\Delta$ in $N$. We say that a map $N_{\QQ} \to \QQ^{k}$ is linear on
a subset $A \subset N_{\QQ}$ if its restriction to $A$ is the
restriction of a linear map.

\begin{definition}
A {\it support map\/} on $\Delta$ is a map $h \colon \vert \Delta \vert
\to \QQ^{k}$ that is linear on every cone $\sigma \in \Delta$.
\end{definition}

For a support map $h \colon \vert \Delta \vert \to \QQ^{k}$, let
$\gamma$ be the cone in $\rq{N} := N \times \ZZ^{k}$ generated by the
graph $\Gamma_{h}$ of $h$, and let $\mathfrak{F}(\gamma)$ denote the
quasifan consisting of all faces of $\gamma$. The {\it filled graph\/} 
of $h$ is the minimal subquasifan $\Lambda_{h}$  of
$\mathfrak{F}(\gamma)$ with $\Gamma_{h} \subset \vert \Lambda_{h}
\vert$. So, $\Lambda_h$ is generated by the cones $\delta \prec
\gamma$ whose relative interior $\delta^{\circ}$ meets $\Gamma_h$.

\begin{definition}\label{convdef}
The support map $h \colon \vert \Delta \vert \to \QQ^{k}$ is called
{\it convex\/}, if the projection $P \colon \rq{N}_{\QQ} \to N_{\QQ}$
is injective on the support $\vert \Lambda_{h} \vert$.
\end{definition}

This notion of convexity includes the classical concept of a convex
support function on a complete fan as defined for example
in~\cite[p.~67]{Fu}:

\begin{remark}\label{suppfunconcomplete}
Let $h \colon \vert \Delta \vert \to \QQ$ be a support map on a fan
$\Delta$. If there are linear forms $u_{\sigma}$, $\sigma \in \Delta$,
on $N$ such that for any pair $\sigma, \tau \in \Delta$ we have
$$ h \vert_{\sigma} = u_{\sigma} \vert_{\sigma}, 
\qquad 
 h \vert_{\tau} \le u_{\sigma} \vert_{\tau}  $$
then $h$ is a convex support map on $\Delta$. Conversely, if $\Delta$
is complete and $h$ is convex then $h$ or $-h$ satisfies the above
condition.
\end{remark}

On noncomplete fans, the concept of convexity for a support
function via the above inequalities is more restrictive than our
concept:

\begin{example}\label{noncompletesuppfct}
Consider the fan $\Delta$ in $\ZZ^2$ generated by
the two maximal cones 
$$\sigma_1:=\cone((1,0),(1,-1)), 
\qquad 
\sigma_2:=\cone((0,1),(1,1))$$
and the support map $h\colon \vert \Delta \vert \to \QQ$ determined
by
$$ h(v_{1},v_{2}) := \left\{
\begin{array}{cl}
2v_{1}+2v_{2} & \quad \text{if } (v_{1},v_{2}) \in \sigma_{1}, \\
-v_{1}+v_{2} &  \quad \text{if } (v_{1},v_{2}) \in \sigma_{2}. \\
\end{array} \right.
$$

Then $h$ is convex: The convex hull $\gamma$ of the graph $\Gamma_{h}$
is a strictly convex cone with four rays, namely  
$$ \gamma = \cone((1,0,2),(1,-1,0),(0,1,1),(1,1,0)). $$
Moreover, the maximal cones of $\Lambda_h$ are precisely the two faces of
$\gamma$ above $\sigma_1$ and $\sigma_2$ respectively. 

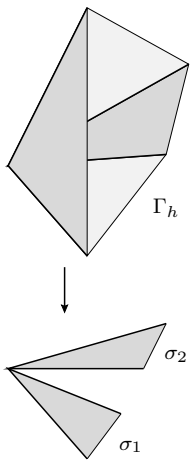
\begin{center}
\begin{picture}(0,0)%
\includegraphics{example_2_4.pstex}%
\end{picture}%
\setlength{\unitlength}{1243sp}%
\begingroup\makeatletter\ifx\SetFigFont\undefined%
\gdef\SetFigFont#1#2#3#4#5{%
  \reset@font\fontsize{#1}{#2pt}%
  \fontfamily{#3}\fontseries{#4}\fontshape{#5}%
  \selectfont}%
\fi\endgroup%
\begin{picture}(3666,9066)(1093,-8644)
\put(4051,-3661){\makebox(0,0)[lb]{\smash{\SetFigFont{8}{9.6}{\familydefault}{\mddefault}{\updefault}
$\Gamma_h$
}}}
\put(3376,-8386){\makebox(0,0)[lb]{\smash{\SetFigFont{8}{9.6}{\familydefault}{\mddefault}{\updefault}
$\sigma_1$
}}}
\put(4276,-6586){\makebox(0,0)[lb]{\smash{\SetFigFont{8}{9.6}{\familydefault}{\mddefault}{\updefault}
$\sigma_2$
}}}
\end{picture}

\bigskip

\end{center}

However neither the function $h$ nor the function $-h$ satisfies the
inequalities of Remark~\ref{suppfunconcomplete}, because we have:
$$ h((0,1)) = 1 < 2, 
\qquad 
   h((1,-1)) = 0 > -2.$$ 
\end{example}

\smallskip

In order to define the notion of strict convexity, we have to note
some observations on convex support maps. The first one is:

\begin{lemma}\label{Sigma_h}
If the support map $h \colon \vert \Delta \vert \to \QQ^{k}$ is
convex, then the projected cones $P(\delta)$, $\delta \in \Lambda_{h}$, form
a quasifan $\Sigma_{h}$ in the lattice $N$.
\end{lemma}

\proof The projection $P$ is injective on any given
$\delta \in \Lambda_h$, and hence induces a bijection between the
faces of $\delta$ and the faces of $P(\delta)$. Moreover,
given $\delta_{1}, \delta_{2} \in \Lambda_{h}$, injectivity of $P$ on
$\vert \Lambda_{h} \vert$ implies 
$$ P(\delta_{1}) \cap P(\delta_{2}) = P(\delta_{1} \cap
\delta_{2}). $$
Since $\delta_{1} \cap \delta_{2}$ is a face of both $\delta_{i}$,
the above consideration yields that $P(\delta_{1}
\cap \delta_{2})$ is a common face of $P(\delta_{1})$ and
$P(\delta_{2})$. \endproof

If $h \colon \vert \Delta \vert \to \QQ^{k}$ is a convex
support map, then
we call $\Sigma_{h}$ the {\it quasifan associated to $h$\/}. 
We need the following properties of this quasifan:

\begin{lemma}\label{Sigmah}
Let $\Sigma_{h}$ be the quasifan associated to a convex support map $h
  \colon \vert \Delta \vert \to \QQ^{k}$. Then we have:
\begin{enumerate}
\item Every cone of $\Delta$ is contained in a cone of $\Sigma_h$.
\item Every cone $\sigma \in \Sigma_{h}$ is generated by the cones
  $\tau \in \Delta$ with $\tau \subset \sigma$. \endproof
\end{enumerate}
\end{lemma}

\begin{definition}
We say that a convex support map  $h \colon \vert \Delta \vert \to
\QQ^{k}$ is {\it strictly convex\/} if its associated quasifan
$\Sigma_h$ equals $\Delta$. 
\end{definition}

Using Remark~\ref{suppfunconcomplete}, one verifies that on a complete
fan $\Delta$, our notion of strict convexity for a support map
$h\colon \vert \Delta \vert \to \QQ$ concides with the usual one, as
defined in~\cite[p.~67]{Fu}. Again, for noncomplete fans the notions
differ:

\begin{example}\label{strconvsuppfct}
The convex support map $h \colon \vert \Delta \vert \to \QQ$ of
Example~\ref{noncompletesuppfct} is even strictly convex.
\end{example}

We now come to the announced main result of this section, namely
the characterization of divisoriality of a toric variety via existence of a
strictly convex support map:

\begin{proposition}\label{tordivchar}
For a fan $\Delta$ in a lattice $N$, the following statements are
equivalent:
\begin{enumerate}
\item $\Delta$ admits a strictly convex support map.
\item The toric variety $X$ associated to $\Delta$ is divisorial.
\end{enumerate}
\end{proposition}

In the proof of this statement, we make use of the following wellknown
characterization of existence of geometric quotients for subtorus
actions in terms of fans, see e.g.~\cite[Theorem 5.1]{Hm}:

\begin{proposition}\label{hamm}
Let $\rq{\Delta}$ be a fan in a lattice $\rq{N}$ with associated toric
variety $\rq{X}$, let $P \colon \rq{N} \to N$ be a surjective lattice
homomorphism, and let $H$ be the subtorus of the big torus of
$\rq{X}$ corresponding to $\ker(P)$. The following statements are
equivalent: 
\begin{enumerate}
\item $P$ is injective on the support $\vert \rq{\Delta} \vert$. 
\item The action of $H$ on $\rq{X}$ has a geometric quotient.
\end{enumerate}
If one of these statements holds, then the quotient variety $\rq{X}/H$ 
is the toric variety determined by the fan $\{P(\sigma); \;
\sigma \in \rq{\Delta}\}$ in $N$. 
\end{proposition} 

\proof[Proof of Proposition~\ref{tordivchar}]
Assume first that the fan $\Delta$ admits a strictly convex
support map $h \colon \vert \Delta \vert \to \QQ^{k}$. 
Then since $\Delta = \Sigma_h$, all cones of $\Sigma_{h}$ are strictly
convex. As before, let $\rq{N} := N \times \ZZ^{k}$. By convexity of
$h$, the projection $P \colon \rq{N}_{\QQ} \to N_{\QQ}$ is an
injection on $\vert\Lambda_h \vert $. In particular, all cones of
$\Lambda_{h}$ are strictly convex. That means that $\Lambda_{h}$ is a
fan. 

The toric variety $\rq{X}$ associated to $\Lambda_{h}$ is
quasiaffine, and the projection $P\colon \rq{N}\to N$ gives 
rise to a toric morphism $p \colon \rq{X} \to X$. According to
Proposition~\ref{hamm},
this toric morphism $p$ is a geometric quotient for the subtorus
action on $\rq{X}$ corresponding to $\ker(P) \subset \rq{N}$.
Thus, Corollary~\ref{divtorvar} yields that $X$ is divisorial.

Suppose now that the toric variety $X$ determined by the fan $\Delta$
is divisorial. By Corollary~\ref{divtorvar}, there is a quasiaffine
toric variety $\rq{X}$ and a toric morphism $p \colon \rq{X} \to X$
such that $H := \ker(p)$ is a subtorus of the big torus of $\rq{X}$
and $p$ is a geometric quotient for the action of $H$ on $\rq{X}$.

Let $p \colon \rq{X} \to X$ arise from a map $P \colon \rq{N} \to N$
of fans $\rq{\Delta}$ and $\Delta$. Since  
$H = \ker(p)$ is  connected, the map $P$ is surjective and
we obtain a section $N \to \rq{N}$ for  $P$. So we may assume that
$\rq{N} = N \times \ZZ^{k}$ holds and that $P$ is the projection onto
the first factor. By the above Proposition~\ref{hamm}, the projection
$P$ is injective on $\vert \rq{\Delta} \vert$. Thus, for each
$\rq{\sigma} \in \rq{\Delta}$, the restriction
$$ P \vert_{\rq{\sigma}} \colon \rq{\sigma} \mapsto \sigma :=
P(\rq\sigma) $$
admits a uniquely determined linear inverse of 
the form $g_{\sigma} = (\id_{N_{\QQ}},
h_{\sigma})$. The maps $h_{\sigma} \colon \sigma \to \QQ^{k}$ patch
together to a support map $h$ on $\Delta$. By construction,  
$\Lambda_{h} = \rq{\Delta}$ and $\Sigma_h=\Delta$.
So $h$ is the desired strictly convex support map on $\Delta$. \endproof

In the remainder of this section we show that convex support maps
in a canonical way define toric morphisms to divisorial toric
varieties. Let $\Delta$ be a fan in a lattice $N$, and let $h \colon
\vert \Delta \vert \to \QQ^{k}$ be a convex support map.

There is a universal method to construct a fan from
the associated quasifan $\Sigma_{h}$:
 Let $\sigma_{\min} \in \Sigma_{h}$ denote its minimal
cone. This is a linear subspace of $N_{\QQ}$. Let $N_{0} :=
\sigma_{\min} \cap N$, set $N_{h} := N/N_{0}$, and denote by $F_{h}
\colon N \to N_{h}$ the projection. The {\it quotient fan\/} of
$\Sigma_{h}$ is the fan
$$ \Delta_{h} := \{F_{h}(\sigma); \; \sigma \in \Sigma_{h} \}. $$
The projection $F_{h} \colon N \to N_{h}$ is a map of the quasifans
$\Sigma_{h}$ and $\Delta_{h}$. Moreover, $F_{h}$ is universal in the
sense that every map of quasifans from $\Sigma_{h}$ to a fan $\Delta'$
factors uniquely through $F_{h}$.

Now, let $X$ and $X_{h}$ denote the toric varieties associated to
the fans $\Delta$ and $\Delta_{h}$ respectively. Our precise statement
is the following:%

\begin{proposition}\label{specialmor} 
The toric variety $X_{h}$ is divisorial, and the projection $F_{h}$
induces a toric morphism $f_{h} \colon X \to X_{h}$.
\end{proposition}

\proof By Lemma~\ref{Sigmah}~i) and the universal property
of the quotient fan $\Delta_h$, the projection $F_{h} \colon N \to
N_{h}$ is a map of the fans $\Delta$ and $\Delta_{h}$ and hence
induces a toric morphism $f_{h} \colon X \to X_{h}$. So we only have
to show that $X_{h}$ is divisorial. In view of
Proposition~\ref{tordivchar}, we look for a strictly convex support
map an $\Delta_{h}$. 

The first step is to construct a strictly convex support map $g$ on
the quasifan $\Sigma_{h}$ associated to $h$: Consider a cone $\sigma
\in \Sigma_{h}$. Then,  as earlier denoting by $P \colon \rq{N} \to N$
the projection, we have $\sigma = P(\delta)$ for some cone $\delta \in
\Lambda_{h}$.

By convexity of $h$, the restriction $P \colon \delta \to \sigma$ has
an inverse of the form $(\id,g_{\sigma})$. The maps $g_{\sigma}$
patch together to a support map $g$ on $\Sigma_h$, and $g$
extends $h$. Moreover, $\Gamma_{g}$ equals $\Lambda_{h}$ and hence the
quasifan associated to $g$ coincides with $\Sigma_h$.

Note that $\Sigma_{g} = \Sigma_{h}$ does not change if
we add a global linear function to $g$.
So we may assume that the support function $g$ vanishes on the minimal
cone of $\Sigma_{g}$. But then we can push down $g$ to a strictly
convex support function on the quotient fan $\Delta_{h}$. \endproof

\section{Toric divisorial reduction}\label{section3}

Fix a toric variety $X$. In~\cite{acha2}, we presented a universal way
to reduce $X$ to a quasiprojective toric variety. In this section we
give an analogous construction, that reduces to divisorial toric
varieties.

\begin{definition}
A {\it toric divisorial reduction\/} of $X$ is a toric morphism $r
\colon X \to X^{\tdr}$ to a divisorial toric variety $X^{\tdr}$ such
that every toric morphism $f \colon X \to Z$ to a divisorial toric
variety $Z$ has a unique factorization $f = \t{f} \circ r$ with a toric
morphism $\t{f} \colon X^{\tdr} \to Z$.
\end{definition}

\begin{theorem}\label{tordivred}
Every toric variety admits a toric divisorial reduction.
\end{theorem}

The proof is given below. We need the following statement on the
pullback of a convex support map:

\begin{lemma}\label{pullback}
Let $F \colon N \to N'$ be a map of fans $\Delta$ and $\Delta'$ in
lattices $N$ and $N'$ respectively. If $h' \colon \vert
\Delta' \vert \to \QQ^{k}$ is a convex support map on $\Delta'$, then 
$h := h' \circ F$ is a convex support map on $\Delta$ and $F$ is a
map of the associated quasifans $\Sigma_{h}$ and $\Sigma_{h'}$.
\end{lemma}

\proof Clearly $h$ is a support map on $\Delta$. To prove convexity of
$h$, we consider the filled graphs $\Lambda_{h}$, $\Lambda_{h'}$ and
the map
$$\rq{F} := F \times \id_{\ZZ^{k}} \colon N \times \ZZ^{k} \to N'
\times \ZZ^{k}. $$
We claim that $\rq{F}$ is a map of the quasifans $\Lambda_{h}$ and
$\Lambda_{h'}$. To verify this, note first that $\rq{F}$ maps the
graph $\Gamma_{h}$ to $\Gamma_{h'}$. Let $\delta \in \Lambda_{h}$. We
have to show that the minimal face $\delta'$ of
$\conv(\Gamma_{h'})$ containing $\rq{F}(\delta)$ belongs to
$\Lambda_{h'}$. Let
$$ G := \id_{\vert \Delta \vert} \times h, \qquad G' :=
\id_{\vert\Delta' \vert} \times h'. $$
By definition of $\Lambda_{h}$, the relative interior $\delta^{\circ}$
of $\delta$ contains a point of the graph of $h$, i.e.
a point of the form $G(v)$ for some 
$v \in \vert \Delta \vert$. By
the choice of $\delta'$ this means $\rq{F}(G(v)) \in (\delta')^{\circ}$. 
On the other hand, by definition of $G$, $G'$ and $\rq{F}$ we have
$$\rq{F}(G(v)) = G'(F(v)) \in \Gamma_{h'}. $$
Hence $\Gamma_{h'} \cap (\delta')^{\circ}\ne\emptyset$.
This implies $\delta' \in\Lambda_{h'}$, and our claim is proved.

For convexity of $h$, we have to show that the
projection $P \colon N \times \ZZ^{k} \to N$ is injective on $\vert
\Lambda_{h} \vert$. Suppose $w_{i} = (v_{i}, t_{i}) \in \vert
\Lambda_{h} \vert$ are two points such that $P(w_{1})$ equals
$P(w_{2})$, that means $v_{1} = v_{2}$. Then we have
$$ P'(\rq{F}(w_{1})) = P'(\rq{F}(w_{2})),$$ 
where $P' \colon N' \times \ZZ^{k} \to N'$ is the projection. Since
$\rq{F}$ is a map of the quasifans $\Lambda_{h}$ and $\Lambda_{h'}$
and $P'$ is injective on $\vert \Lambda_{h'} \vert$, this
implies $\rq{F}(w_{1}) = \rq{F}(w_{2})$. In particular, we have $t_{1}
= t_{2}$ and thus $w_{1} = w_{2}$.

Finally, the fact that $F$ is a map of the quasifans $\Sigma_{h'}$ and
$\Sigma_{h}$ follows immediately from the fact that $\rq{F}$ is a map
of the quasifans $\Lambda_{h'}$ and $\Lambda_{h}$. \endproof

\proof[Proof of Theorem~\ref{tordivred}]
Let $X$ be a toric variety arising from a fan $\Delta$ in a
lattice $N$. First we show that any given toric morphism $f\colon X
\to Z$ from $X$ to a divisorial variety $Z$ factors uniquely through one
of the toric morphisms $f_h$ arising from  a convex support map on $\Delta$
as in Proposition~\ref{specialmor}.

To see this, consider the map of fans $F\colon \Delta\to\Delta'$
associated to the given toric morphism $f$ and 
choose a strictly convex support map $h'$ on 
$\Delta'$. Lemma~\ref{pullback} tells us that by pulling back $h'$ via  $F$, 
we obtain a convex support map $h$  on $\Delta$.
Moreover, $F$ defines a map of quasifans from 
$\Sigma_h$ to $\Sigma_{h'} = \Delta'$. 

Now, the map of fans $F$ factors as a map of fans through the
projection $F_{h} \colon N \to N_{h}$, i.e., $F$ induces a map from
the quotient fan $\Delta_{h}$ of $\Sigma_{h}$ to $\Delta'$. Obviously,
the corresponding toric morphism is the desired factorization of $f \colon
X \to Z$
 through $f_{h} \colon X \to X_{h}$.

Now let us take a closer look at the toric morphisms $f_h \colon X \to
X_{h}$ arising from convex support maps. Recall that 
the morphism $f_h$ is already determined by the quasifan $\Sigma_h$
associated to $h$. By Lemma~\ref{Sigmah}~ii), each such quasifan has
the property that all cones are generated by cones of $\Delta$.
Consequently there exist only finitely many
 of such quasifans, say
$\Sigma_1,\dots,\Sigma_r$.

Let $f_i \colon X \to Y_i$ denote the toric morphisms
to divisorial toric varieties determined by $\Sigma_i$,
and consider their product $f:=f_1\times\ldots\times f_r$.
Let $Y$ denote the closure of the image $f(X)$ in
$Y_1\times\dots\times Y_r$. The normalization $\t{Y}$ of $Y$ is again a 
divisorial toric variety, and $f$ lifts to a toric morphism to $\t{Y}$.
In $\t{Y}$ we choose the smallest open
toric subvariety $Y'$ containing the image of $f$, and restricting $f$,
we obtain a toric morphism $r\colon X\to Y'$. 

By construction, for every $i$ we have a unique factorization of $f_i$
through $r$, namely $f_i=\pr_i\circ r$, where $\pr_i\colon Y'\to Y_i$
denotes the restriction of the projection on the $i$-th factor.
This proves that $r$ is the desired toric divisorial reduction.  
\endproof

We conclude this section with some examples. Note that
any two-dimensional toric variety is simplicial and hence
divisorial. So the minimal dimension for interesting
examples is $3$.

\begin{example}\label{beispiel1}
If a toric variety does not admit nontrivial
effective Cartier divisors, see e.g.~\cite[p.~25]{Fu}, then its toric
divisorial reduction is a point. 
\end{example}

\begin{example}\label{beispiel2}%
Consider the following eight vectors in $\QQ^3$:
$$
\begin{array}{lclclcl}
v_1 := (2,2,1), & & v_2 := (-2,2,1), & & v_3 := (-2,-2,1),
& & v_4:=(2,-2,1), \\
v_5:=(1,1,1), & & v_6:=(-1,1,1), & & v_7:=(-1,-1,1), 
& & v_{8} := (2/3,1/3,1).
\end{array}
$$ 
Let $\Delta$ denote the fan in $\ZZ^{3}$ with maximal cones 
$$
\begin{array}{lcl}
\sigma_1:=\cone(v_1,v_2,v_5,v_6), & & \sigma_2:=\cone(v_2,v_3,v_6,v_7),\\
\sigma_3:=\cone(v_3,v_4,v_7,v_8), & & \sigma_4:=\cone(v_1,v_4,v_5,v_8),\\
\sigma_5:=\cone(v_5,v_6,v_7,v_8). & & \\
\end{array} 
$$

\goodbreak

\begin{center}
\vbox{\begin{picture}(0,0)%
\includegraphics{example_3_5.pstex}%
\end{picture}%
\setlength{\unitlength}{1184sp}%
\begingroup\makeatletter\ifx\SetFigFont\undefined%
\gdef\SetFigFont#1#2#3#4#5{%
  \reset@font\fontsize{#1}{#2pt}%
  \fontfamily{#3}\fontseries{#4}\fontshape{#5}%
  \selectfont}%
\fi\endgroup%
\begin{picture}(9045,3000)(721,-3191)
\put(7201,-3161){\makebox(0,0)[lb]{\smash{\SetFigFont{7}{8.4}{\familydefault}{\mddefault}{\updefault}
$v_4$
}}}
\put(3196,-491){\makebox(0,0)[lb]{\smash{\SetFigFont{7}{8.4}{\familydefault}{\mddefault}{\updefault}
$v_2$
}}}
\put(721,-3191){\makebox(0,0)[lb]{\smash{\SetFigFont{7}{8.4}{\familydefault}{\mddefault}{\updefault}
$v_3$
}}}
\put(4216,-2536){\makebox(0,0)[lb]{\smash{\SetFigFont{7}{8.4}{\familydefault}{\mddefault}{\updefault}
$\sigma_3$
}}}
\put(7456,-1876){\makebox(0,0)[lb]{\smash{\SetFigFont{7}{8.4}{\familydefault}{\mddefault}{\updefault}
$\sigma_4$
}}}
\put(7561,-1106){\makebox(0,0)[lb]{\smash{\SetFigFont{7}{8.4}{\familydefault}{\mddefault}{\updefault}
$v_5$
}}}
\put(4021,-1091){\makebox(0,0)[lb]{\smash{\SetFigFont{7}{8.4}{\familydefault}{\mddefault}{\updefault}
$v_6$
}}}
\put(5791,-1036){\makebox(0,0)[lb]{\smash{\SetFigFont{7}{8.4}{\familydefault}{\mddefault}{\updefault}
$\sigma_1$
}}}
\put(5101,-1741){\makebox(0,0)[lb]{\smash{\SetFigFont{7}{8.4}{\familydefault}{\mddefault}{\updefault}
$\sigma_5$
}}}
\put(9766,-626){\makebox(0,0)[lb]{\smash{\SetFigFont{7}{8.4}{\familydefault}{\mddefault}{\updefault}
$v_1$
}}}
\put(6466,-2276){\makebox(0,0)[lb]{\smash{\SetFigFont{7}{8.4}{\familydefault}{\mddefault}{\updefault}
$v_8$
}}}
\put(2731,-2531){\makebox(0,0)[lb]{\smash{\SetFigFont{7}{8.4}{\familydefault}{\mddefault}{\updefault}
$v_7$
}}}
\put(2746,-1756){\makebox(0,0)[lb]{\smash{\SetFigFont{7}{8.4}{\familydefault}{\mddefault}{\updefault}
$\sigma_2$
}}}
\end{picture}

\bigskip

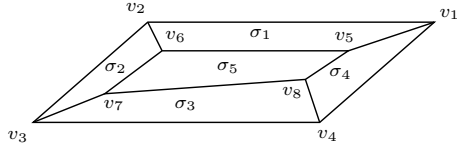
{\small Intersection of $\Delta$ with the plane $x_{3} = 1$.}}
\end{center}

The identity on $\ZZ^3$ defines a map of fans from $\Delta$
to the fan of faces $\mathfrak{F}(\sigma)$ of the cone
$\sigma:=\cone(v_1,v_2,v_3,v_4)$. We claim that the corresponding
toric morphism $r\colon X_{\Delta}\to X_{\sigma}$  is the toric
divisorial reduction of $X_{\Delta}$. 

To see this, consider a convex support map $h\colon|\Delta|\to\QQ^k$,
and its associated quasifan $\Sigma_{h}$. Lemma~\ref{Sigmah}
implies that we have only two possibilities, namely
$\Sigma_h=\mathfrak{F}(\sigma)$ or $\Sigma_h=\Delta$. 
Thus, to verify our claim, we only have to exclude the latter
possibility, i.e., we have to show that $h$ cannot be strictly
convex.

Otherwise, let $\delta_{5} \in \Lambda_{h}$ be the maximal
cone above $\sigma_{5}$ and choose a linear form
$\lambda\colon N_{\QQ}\times\QQ^k\to\QQ$ that is nonnegative
on $\gamma := \conv(\Gamma_{h})$ and fulfills $\delta_5 = \gamma \cap
\lambda^{\perp}$. Pulling back $\lambda$ via $\id_N\times h$, we obtain
a nonnegative support function $g$ on $\Delta$ vanishing
precisely on $\sigma_5$. Note that
$$ g(v_{1}) = g(v_{2}) = g(v_{3}) .$$
Moreover, we have the relations
$$ v_{4} = 17v_{3} -28v_{7} + 12v_{8}, \qquad v_{4} = 5v_{1} 
-16v_{5} + 12v_{8}. $$
Applying $g$, we obtain $17g(v_{3}) = 5g(v_{1})$. This contradicts
$g(v_{1}) = g(v_{3})$. So, $h$ cannot be
strictly convex and our claim is proved.
\end{example}

\begin{example}\label{nichtsurjtdr}%
We describe a toric variety with a nonsurjective toric divisorial
reduction. Similarly to the preceding example, consider the vectors
$$
\begin{array}{lclcl}
v_1 := (2,2,1,0), & & v_2 := (-2,2,1,0), & & v_3 := (-2,-2,1,0), \\
v_4:=(2,-2,1,0), & & v_5:=(1,1,1,0), & & v_6:=(-1,1,1,0), \\
v_7:=(-1,-1,1,0), & & v_{8} := (2/3,1/3,1,0). & &  \\
\end{array}
$$ 
in $\QQ^{4}$ and let furthermore $e_{4}$ be the fourth canonical 
base vector. Let
$\Delta$ denote the fan in $\ZZ^{4}$ with maximal cones
$$
\begin{array}{lcl}
\sigma_1:=\cone(v_1,v_2,v_5,v_6), & & \sigma_2:=\cone(v_2,v_3,v_6,v_7),\\
\sigma_3:=\cone(v_3,v_4,v_7,v_8), & & \sigma_4:=\cone(v_1,v_4,v_5,v_8),\\
\sigma_5:=\cone(v_5,v_6,v_7,v_8). & & \sigma_{6} := \cone (v_{5},v_{6},e_{4})\\
\end{array} 
$$

\begin{center}
\vbox{\begin{picture}(0,0)%
\includegraphics{example_3_6.pstex}%
\end{picture}%
\setlength{\unitlength}{1184sp}%
\begingroup\makeatletter\ifx\SetFigFont\undefined%
\gdef\SetFigFont#1#2#3#4#5{%
  \reset@font\fontsize{#1}{#2pt}%
  \fontfamily{#3}\fontseries{#4}\fontshape{#5}%
  \selectfont}%
\fi\endgroup%
\begin{picture}(9240,5433)(571,-5161)
\put(7201,-5161){\makebox(0,0)[lb]{\smash{\SetFigFont{7}{8.4}{\familydefault}{\mddefault}{\updefault}
$v_4$
}}}
\put(571,-5146){\makebox(0,0)[lb]{\smash{\SetFigFont{7}{8.4}{\familydefault}{\mddefault}{\updefault}
$v_3$
}}}
\put(2701,-4561){\makebox(0,0)[lb]{\smash{\SetFigFont{7}{8.4}{\familydefault}{\mddefault}{\updefault}
$v_7$
}}}
\put(6256,-4306){\makebox(0,0)[lb]{\smash{\SetFigFont{7}{8.4}{\familydefault}{\mddefault}{\updefault}
$v_8$
}}}
\put(7831,-3616){\makebox(0,0)[lb]{\smash{\SetFigFont{7}{8.4}{\familydefault}{\mddefault}{\updefault}
$v_5$
}}}
\put(3196,-3421){\makebox(0,0)[lb]{\smash{\SetFigFont{7}{8.4}{\familydefault}{\mddefault}{\updefault}
$v_6$
}}}
\put(9811,-2731){\makebox(0,0)[lb]{\smash{\SetFigFont{7}{8.4}{\familydefault}{\mddefault}{\updefault}
$v_1$
}}}
\put(2821,-2671){\makebox(0,0)[lb]{\smash{\SetFigFont{7}{8.4}{\familydefault}{\mddefault}{\updefault}
$v_2$
}}}
\end{picture}

\bigskip

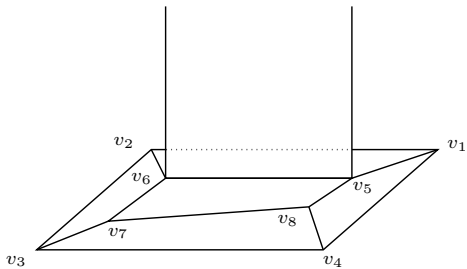
{\small Intersection of $\Delta$ with the hyperplane $x_{3} = 1$.}}
\end{center}

The identity on $\ZZ^4$ defines a map of fans from $\Delta$
to the fan of faces $\mathfrak{F}(\sigma)$ of the cone
$\sigma:=\cone(v_1,v_2,v_3,v_4, e_{4})$. We claim that the corresponding
toric morphism $r\colon X_{\Delta}\to X_{\sigma}$  is the toric 
divisorial reduction of $X_{\Delta}$. Note that this map is not
surjective.

Let us verify the claim. If $h$ is a convex support map it follows that
$|\Sigma_h|\subset \sigma$. The restriction of $h$ to the support
of the subfan $\Delta'$ of $\Delta$ generated by the cones
$\sigma_1,\ldots,\sigma_5$
defines a convex support map $h'$ of $\Delta'$.
So by the previous example, $\Sigma_{h'}=\mathfrak{F}(\sigma')$,
where $\sigma'$ denotes the cone generated by $v_1,\dots,v_4$.

Now Lemma~\ref{pullback} implies that the smallest cone 
$\tau$ in $\Sigma_h$ containing $\sigma_5$
also contains all of $\sigma'$. That means by Lemma~\ref{Sigmah} that
either $\tau=\sigma'$ or $\tau=\sigma$. In any case, since
 $\sigma'$ is a face of $\sigma$ we obtain $\sigma'\in \Sigma_h$. 

Next consider the smallest cone $\tau'\in\Sigma_h$ containing $\sigma_6$.
We have $v_5,v_6\in \sigma_6$, so the cone $\tau'$ meets $\sigma'$ in its
relative interior. Since $\Sigma_h$ is a quasifan, we can conclude that
$\sigma'$ is in fact a face of $\tau'$. Because $e_4\in\tau$ this implies
$\tau'=\sigma$, and we obtain $\Sigma_h=\mathfrak{F}(\sigma)$.
\end{example}

\section{A Lifting Lemma}\label{section4}

Here we relate regular maps between divisorial toric
prevarieties to regular maps between quasiaffine toric varieties.
For maps of projective spaces, this is a classical observation:

\begin{example}
Let $f \colon \PP_{n} \to \PP_{m}$ be a regular map of projective
spaces. Then $f$ is of the form
$$ [z_{0}, \ldots, z_{n}] \mapsto 
[f_{0}(z_{0}, \ldots, z_{n}), \ldots, f_{m}(z_{0}, \ldots, z_{n})]$$
with homogeneous polynomials $f_{i}$ that are pairwise of the same
degree. In other words, there is a lifting
$$\xymatrix{
{\KK^{n+1} \setminus \{0\}} \ar[r]^{\rq{f}} \ar[d] & 
{\KK^{m+1} \setminus \{0\}} \ar[d] \\
{\PP_{n}} \ar[r]^{f} & {\PP_{m}}
}$$
\end{example}

The main result of this section is the following generalization of the
above lifting statement:

\begin{lemma}\label{lifting}
Let $f\colon X_{1} \to X_{2}$ be a regular map of divisorial toric
prevarieties such that $f(X_{1})$ intersects the big torus of $X_{2}$.
Then there exists a commutative diagram
$$ \xymatrix{
{\h{X}_{1}} \ar[r]^{\h{f}} \ar[d]_{q_{1}} &
{\h{X}_{2}} \ar[d]^{q_{2}} \\
{X_{1}} \ar[r]^{f} & X_{2}
}$$
where  $\h{X}_{1}$, $\h{X}_{2}$ are quasiaffine toric varieties,
$q_{i} \colon \h{X}_{i} \to X_{i}$ are geometric prequotients for free
subtorus actions on $\h{X}_{i}$ and $\h{f} \colon \h{X}_{1} \to \h{X}_{2}$
is a regular map.
\end{lemma}

\proof We use the ideas and methods presented
in~\cite[Section~2]{ha}. Choose effective $T_{i}$-invariant Cartier
divisors $D^{i}_{1},  \ldots, D^{i}_{r_{i}}$ on $X_{i}$ such that the
complements $X_{i} \setminus \Supp(D^{i}_{j})$ form an affine cover of
$X_{i}$. Let $W_{i} \subset \CDiv(X_{i})$ denote the subgroup
generated by $D^{i}_{1}, \ldots, D^{i}_{r_{i}}$. The pullback via $f$
gives rise to a group homomorphism
$$ \psi \colon W_{2} \to \CDiv(X_{1}), \qquad D \mapsto f^{*}(D). $$

Enlarge $W_{1}$ by adding the image $\psi (W_{2})$. Note that the 
line bundles determined by the divisors of $W_{i}$ are
$T_{i}$-linearizable, see~\cite[p.~67, Remark]{DMV}. We shall regard
$\psi$ in the sequel as a homomorphism from $W_{2}$ to
$W_{1}$. Consider the $\mathcal{O}_{X_{i}}$-algebras
$$ \mathcal{A}_{i} := \bigoplus_{D \in W_{i}} \mathcal{O}_{D}(X_{i}) $$ 
and their associated relative spectra $\h{X}_{i} :=
\Spec(\mathcal{A}_{i})$. By \cite[Remark~2.1]{ha}, the inclusion
$\mathcal{O}_{X_{i}} \subset \mathcal{A}_{i}$ gives rise to a
geometric prequotient $q_{i} \colon \rq{X}_{i} \to X_{i}$ for the
free action of the algebraic torus $H_{i} := \Spec(\KK[W_{i}])$ on
$\h{X}_{i}$ induced by the $W_{i}$-grading of $\mathcal{A}_{i}$. 

Since $W_{1}$ and $W_{2}$ define ample groups of line bundles in the
sense of \cite[Definition~2.2]{ha}, each $\h{X}_{i}$ is in fact a
quasiaffine variety. Moreover, by \cite[Proposition~2.3]{ha}, the
variety $\rq{X}_{i}$ carries a regular action of the algebraic torus
$T_{i}$ commuting with the action of $H_{i}$ such that $q_{i} \colon
\h{X}_{i} \to X_{i}$ becomes $T_{i}$-equivariant. It follows that
$\h{X}_{i}$ is a toric variety with big torus $\h{T}_{i} = T_{i}
\times H_{i}$.

We still have to construct the lifting $\h{f} \colon \h{X}_{1} \to
\h{X}_{2}$. As to this, note that for every affine open subset $U
\subset X_{2}$, we obtain a homomorphism of $W_{i}$-graded algebras by
setting 
$$ \mathcal{A}_{2}(U) \to \mathcal{A}_{1}(f^{-1}(U)), \qquad
\mathcal{O}_{D}(U) \ni h \mapsto f^{*}(h) \in
\mathcal{O}_{\psi(D)}(U) \qquad (D\in W_2). $$

Note that on the homogeneous component $\mathcal{A}_{2}(U)_{0}$, this
is just the comorphism of the map $f$. By definition of $\h{X}_{i}$
and the maps $q_{i} \colon \h{X}_{i} \to X_{i}$, each of the above
homomorphisms gives rise to a lifting
$$ \h{f}_{U} \colon q_{1}^{-1}(f^{-1}(U)) \to q_{2}^{-1}(U)$$
of the restriction $f \colon f^{-1}(U) \to U$. By construction, the
maps $\h{f}_{U}$ patch together to the desired lifting $\h{f} \colon
\rq{X}_{1} \to X_{2}$ of $f \colon X_{1} \to X_{2}$. \endproof

The following observation will be needed later to obtain equivariance
properties for the lifting $\h{f} \colon \h{X}_{1} \to \h{X}_{2}$
constructed in the above Lemma.

\begin{lemma}\label{equivariance}
For $i=1,2$, let $T_{i}$ be algebraic tori and let $Y_{i}$ be
irreducible $T_{i}$-varieties such that $T_{2}$ acts freely on
$Y_{2}$. If $f \colon Y_{1} \to Y_{2}$ is regular and maps the orbits
of $T_{1}$ into orbits of $T_{2}$, then there is a homomorphism
$\varphi \colon T_{1} \to T_{2}$ such that $f(t \mal x) = \varphi(t)
\mal f(x)$ holds for all $(t,x) \in T_{1} \times Y_{1}$.
\end{lemma}

\proof By Sumihiro's Theorem~\cite[Corollary 2]{Su}, we may assume
that $Y_{2}$ is affine. Thus, there is an algebraic quotient $Y_{2}
\to Y$ for the action of $T_{2}$ on $Y_{2}$. 
Since $T_{2}$ acts freely, the quotient map $Y_{2} \to Y$ is
equivariantly locally trivial. Thus, shrinking $Y$, we may even assume
that $Y_{2} = T_{2} \times Y$ holds. In particular, one has $f =
(f_{1}, f_{2})$ with regular maps $f_{1} \colon Y_{1} \to T_{2}$ and
$f_{2} \colon Y_{1} \to Y$. So, we obtain a regular map
$$ \Phi \colon T_{1} \times Y_{1} \to T_{2}, \qquad (t,x) \mapsto
f_{1}(t \mal x) f_{1}(x)^{-1}.   $$ 

For fixed $x \in Y_{1}$, the map $t \mapsto \Phi(t,x)$ maps
the neutral element of $T_{1}$ to the neutral element of $T_{2}$ and
hence is necessarily a homomorphism of the tori $T_{1}$ and
$T_{2}$. By rigidity of tori \cite[III.8.10]{Borel}, the map $\Phi$ does
not depend on $x$. So there is a homomorphism $\varphi \colon T_{1}
\to T_{2}$ with $\Phi(t,x) = \varphi(t)$ for all $(t,x) \in T_{1}
\times Y_{1}$. Clearly, $\varphi$ is as desired. \endproof

A  different aspect of the lifting problem is discussed extensively in
\cite{Be}: Given two quotient presentations $\h{X}_{i} \to X_{i}$ of
toric varieties in the sense of \cite{achasc} and a regular map $f
\colon X_{1} \to X_{2}$,  when can  this map be lifted to a
regular map $F \colon \h{X}_{1} \to \h{X}_{2}$?

\section{Decomposition of regular maps}\label{section5}

Let $X$ be a toric variety with big torus $T$ and 
consider the action of a closed subgroup $H \subset T$ on $X$. Here we
provide  the key to relate $H$-invariant regular maps $X \to Y$ to
$H$-invariant toric morphisms:

\begin{lemma}\label{Zerlegungssatz}
Let $f \colon X \to Y$ be an $H$-invariant regular map to a
  divisorial variety~$Y$. Then there exists a dominant
  $H$-invariant toric morphism $g \colon X \to X'$ to a divisorial
  toric variety $X'$, an open subset $U \subset X'$ with 
  $g(X) \subset U$ and a regular map $h \colon U \to
  Y$ such that $f = h \circ g$.
\end{lemma}

\proof First we reduce the problem to the case that $H$ is
connected. Suppose that $g \colon X \to X'$ and $h \colon U \to Y$
satisfy the assertion for the identity component $H^{0}$ of $H$. Then
$g$ induces an action of the finite abelian group $\Gamma := H/H^{0}$
on $X'$. Let $p \colon X' \to X''$ be the geometric quotient for this
action. Note that $p$ is a toric morphism. Using
Corollary~\ref{divtorvar}, we see that the variety $X''$ is again
divisorial.

By appropriate shrinking, we achieve that $U$ is
$\Gamma$-invariant. Since $p$ is geometric, $p(U)$ is open in $X''$
and the restriction $p \colon U\to p(U)$ is again a geometric quotient
for the action of $\Gamma$. Since $h$ is $\Gamma$-invariant, we have
$h = h' \circ p$ for some regular map $h' \colon p(U) \to Y$. It
follows that $f = h' \circ (p \circ g)$ is the desired
decomposition. Consequently, it suffices to give the proof for
connected $H$.

The next simplification provides the link to the toric setting:
As mentioned before, we can realize $Y$ as a closed subvariety of 
a smooth toric prevariety $Z$ of affine intersection,
see~\ref{einbettungssatz}.
Let $Z' \subset Z$ denote the minimal orbit closure of
the big torus of $Z$ such that $f(X) \subset Z'$ holds. Then $Z'$
is again a smooth toric
prevariety of affine intersection, but in $Z'$
the image $f(X)$ intersects the big torus.

Now,  for the moment regard $f$ as a map from $X$ to $Z'$ and suppose
that $g \colon X \to X'$ and $h \colon U \to Z'$ satisfy the assertion
for $f \colon X \to Z'$. Taking closures in $U$ and $Z'$
respectively, we obtain
$$ h(U) \subset h\left(\b{g(X)}\right) \subset \b{h(g(X))} = \b{f(X)}
\subset Y. $$

That means $h$ is in fact a map from $U$ to $Y$. Thus $X'$, $g$, $h$
and $U$ also provide the desired data for the original $f \colon X \to
Y$. Consequently, we can assume in the sequel that $Y$ is a smooth
toric prevariety of affine intersection and that $f(X)$ intersects the
big torus of $Y$. But then according to Lemma~\ref{lifting} there
is a commutative diagram
$$\xymatrix{
{\h{X}} \ar[r]^{\h{f}} \ar[d]_{p} &
{\h{Y}} \ar[d]^{q} \\
{X} \ar[r]^{f} & Y}$$
where $\h{X}$, $\h{Y}$ are quasiaffine toric varieties and the
vertical maps are geometric prequotients for free actions of subtori
$H_{X}$ and $H_{Y}$ of the big tori of $\h{X}$ and $\h{Y}$
respectively. We may even assume that $\h{X} =X$ holds:

Let $H' :=
p^{-1}(H)$ and suppose that the $H'$-invariant regular map $f' := f
\circ p$ admits a decomposition of the form $f' = h'
\circ g'$ with a dominant $H'$-invariant toric morphism $g' \colon
\h{X} \to X'$ and a regular map $h' \colon U \to Y$ defined on an open
neighbourhood $U$ of the image of $g'$.

Then, by the universal property of $p$, there is a toric morphism $g
\colon X \to X'$ with $g' = g \circ p$. Clearly this morphism is
dominant. Moreover, since $p$ is surjective, it is $H$-invariant and
$g(X) \subset U$ holds. Consequently, $f = h' \circ g$ is a
decomposition as wanted. So it suffices to prove the assertion for the
case that $\h{X} = X$ and $H_X=1$ hold  and $p$ is the identity map.

Now we consider the regular map $\h{f}\colon X\to \h{Y}$ as
a map from an $H$-variety to an $H_Y$-variety. Since $q\circ\h{f}=f$
is $H$-invariant, every  $H$-orbit is mapped by $\h{f}$ into
a fiber of $q$. On the other hand, the fibers of $q$ are precisely
the $H_Y$-orbits. So we can apply Lemma~\ref{equivariance} and
conclude that $\h{f}$ is $H$-equivariant with respect to a
homomorphism $H\to H_Y$.

Choosing a locally closed toric embedding $\h{Y} \subset
\CC^{s}$, we obtain a homomorphism $H_Y\to\CC^s$,
and the induced map $\h{f}\colon X \to \CC^{s}$ is $H$-equivariant
with respect to the homomorphism $H\to H_Y\to\CC^s$. So the components of
$\h{f}$ are $H$-homogeneous regular functions.
By writing the components of $\h{f}$ as linear combinations of
character functions of the big torus $T \subset X$, and using the 
summands to define a toric morphism $g' \colon X \to \CC^r$,
we obtain a
decomposition of $\h{f}$ in the form $\h{f} = s \circ g'$, with 
 a linear map $s \colon \CC^r \to \CC^{s}$.
Note that $g'$ induces an action of $H$ on $\CC^{r}$ making $s
\colon \CC^{r} \to \CC^{s}$ into an $H$-equivariant map.

Let $W$ be the normalization of the closure of $g'(X)$ in $\CC^r$.
Then $W$ is an affine toric variety with big torus $g'(T)$. We
can lift
$g'$ to a dominant toric morphism $\h{g} \colon
X \to W$, and pull back $s$ to a regular map $\h{s} \colon W \to
\CC^s$. Both, $\h{g}$ and $\h{s}$, are again equivariant for the
induced $H$-action on $W$. The set $V := \h{s}^{-1}(\h{Y})$ is
$H$-invariant and open in $W$. Moreover, we have $\h{g}(X) \subset V$.
So far, we are in the following situation:
$$\xymatrix{
{V} \ar[r]^{\h{s}} &
{\h{Y}} \ar[d]^{q} \\
{X} \ar[u]^{\h{g}} \ar[r]^{f} & 
{Y} }$$

Since $\h{s} \colon W \to \CC^{s}$ is an affine map, also its
restriction $\h{s} \colon V \to \h{Y}$ is affine. Thus $q \circ
\h{s} \colon  V \to Y$ is an affine $H$-invariant regular map. 
Existence of an affine $H$-invariant map $V \to Y$ already implies
existence of a good quotient $p \colon V \to V \quot H$ for the action
of $H$, see e.g.~\cite[Prop. 3.12]{Ne}. So we obtain the following
commutative diagram of regular maps: 
$$\xymatrix{
{V} \ar[r]^{p} & 
{ V\quot H} \ar[d]^{h} \\
{X} \ar[u]^{\h{g}} \ar[r]^{f} &
{Y} }$$

Note that $g:=p\circ\h{g}\colon X\to V\quot H$ is $H$-invariant
and $V\quot H$ is divisorial, because $Y$ is divisorial
and $h$ is an affine morphism. So the decomposition $f=h\circ g$
is almost as wanted. To complete the proof it suffices to show
that we can embed $V\quot H$ as an open subset into a divisorial
toric variety $X'$ such that $g$ viewed as a morphism from 
$X$ to $X'$ is toric.

For this last step we argue as follows: Note that we constructed $V$
as an open $H$-invariant subset of 
the toric variety $W$. In \cite{Sw}, J.~\'Swi\c{e}cicka shows that
``maximal'' open subsets with a good quotient by a given subtorus
in a toric variety are in fact toric subvarieties.

More precisely, according to \cite[Corollary~2.4]{Sw}, $V$  is
contained in an open toric subvariety $V' \subset W$ with a good toric
quotient $p' \colon V' \to V' \quot H$ such that  the induced map 
$V \quot H \to V' \quot H$ is an open
inclusion. Of course, we can choose $V'$ in such a manner that $V'
\quot H = T' \mal (V \quot H)$ holds, where $T'$ denotes the big torus
of $V' \quot H$. We set $X':=V'\quot H$ and $U:=V\quot H$ and arrive at
the following commutative diagram: 
$$\xymatrix{
 & {V'} \ar[r]^{p'} 
 & {V'\quot H} \ar@{}[r]|=
 & X' \\
{X} \ar[r]^{\h{g}} 
 & {V} \ar@{}[u]|\cup \ar[r]^{p} 
 & {V \quot H}  \ar@{}[u]|\cup \ar@{}[r]|=
 & {U} \ar@{}[u]|\cup
}$$

The morphism $X \to V'$ sending $x$ to $\rq{g}(x)$ is a dominant toric
morphism because $\rq{g} \colon X \to W$ is one. Hence the same is
true for $g = p'\circ \h{g} \colon X \to X'$. Moreover, because
$\h{g}(X) \subset V$ holds, we conclude that the big torus $T'$ of
$X'$ is contained in $U$. It follows that the complement 
$X' \setminus U$ is of codimension at least 2 in $X'$. Thus
Lemma~\ref{translates} yields that the toric variety $X'$ is
also divisorial. This ends the proof. \endproof

\section{Divisorial reduction and categorical quotients}\label{section6}

In this section we come to the main results of this article. Recall
from \cite{Mu} that a {\it categorical quotient\/} for a 
$G$-variety $X$ is a $G$-invariant regular map $X \to Y$ such that any
$G$-invariant regular map $X \to Z$ factors uniquely through $X \to
Y$. Clearly, this notion can be restricted to any subcategory of the
category of algebraic varieties, as soon as the $G$-variety $X$
belongs to this subcategory.

We give an answer to the problem of existence of categorical quotients
for subtorus actions in the divisorial category. Our method of proof in
fact solves the existence problem of a more general universal object:
Consider a toric variety $X$ with big torus $T$ and the action of a
subtorus $H \subset T$.

\begin{definition}\label{divreddef}
An {\it $H$-invariant divisorial reduction\/} of $X$ is a regular
map $r \colon X \to Y$ to a divisorial variety
$Y$ such that every $H$-invariant regular map $f \colon X \to Z$ 
to a divisorial variety $Z$ admits a unique factorization $f
= \t{f} \circ r$ with a regular map $\t{f} \colon Y \to Z$.
If $H=1$, then we simply speak of a {\it divisorial reduction\/}.
\end{definition}

A candidate for such a reduction is constructed in two steps.
First, recall from~\cite{acha1} that there is a toric quotient for the
action of $H$ on $X$, that means a toric morphism
$$  p \colon X \to X \tq H $$
which is a categorical quotient for the action of $H$ on $X$ in the
category of toric varieties. In a second step, construct the toric
divisorial reduction of the toric quotient space as described in
Section~\ref{section3}:
$$ q \colon X \tq H \to (X \tq H)^{\tdr} .$$

\begin{theorem}\label{Hinvardivred}
For a toric variety $X$, the following statements are equivalent:
\begin{enumerate}
\item $X$ admits an $H$-invariant divisorial reduction.
\item The composition $q \circ p \colon X \to Z$ is surjective.
\end{enumerate}
Moreover, if one of these statements holds, then $q \circ p$ is
the $H$-invariant divisorial reduction.
\end{theorem}

Applying this result to divisorial toric varieties $X$, we obtain the
following solution for the above quotient problem:

\begin{corollary}\label{divcatquot}
The action of a subtorus $H$ on a divisorial toric 
variety $X$ admits a categorical quotient in the category of divisorial
varieties if and only if the composition of $X \to X \tq H$ and $X \tq H
\to (X \tq H)^{\tdr}$ is a surjective map. \endproof
\end{corollary}

A further special case of Theorem~\ref{Hinvardivred} is the case of a
trivial torus $H=1$. Here we obtain the following: 

\begin{corollary}\label{divredchar}
A toric variety admits a divisorial reduction if and only if its toric
divisorial reduction is surjective. \endproof 
\end{corollary}

\proof[Proof of Theorem~\ref{Hinvardivred}]
Assume first that $q \circ p$ is surjective. We show that a
given $H$-invariant regular map $f \colon X \to Z$ to a divisorial
variety $Z$ factors through $q \circ p$. Lemma~\ref{Zerlegungssatz}
yields a decomposition $f = h \circ g$ with an $H$-invariant dominant
toric morphism $g \colon X \to X'$ to a divisorial toric variety
$X'$.

By the universal properties of $p$ and $q$, the toric morphism $g$ has
a factorization $g = g' \circ (q \circ p)$. By
surjectivity of $q \circ p$, the map $h$ is defined on a neighbourhood
of the image of $g'$. Hence $f = (h \circ g') \circ (q \circ p)$ is
the desired factorization. Thus $q \circ p$ is the $H$-invariant
divisorial reduction of $X$.

Conversely, suppose that $X$ has an $H$-invariant divisorial reduction
$r \colon X  \to Y$. Since the normalization of a divisorial variety is
again divisorial, we can conclude that $Y$ is normal. Moreover,
 the universal property of $r \colon X \to Y$ implies that $r$
 is surjective, and that $Y$
 inherits a set-theoretical action of the big torus $T \subset X$
 making $r$ equivariant. Note that a priori it is not clear that this
 action is regular, so we cannot treat $Y$ as a toric variety.

Let $Z := (X \tq H)^{\tdr}$. We  shall compare the $H$-invariant
divisorial reduction $r \colon X \to Y$ with the toric morphism $q
\circ p \colon X \to Z$. On the one hand, because of the universal
property of $r$, the map $q\circ p$ factors uniquely through $r$. So
there is a unique regular map $\alpha \colon Y \to Z$ with $q \circ
p = \alpha \circ r$.

On the other hand,
Lemma~\ref{Zerlegungssatz} provides a decomposition $r = h \circ g$
with a dominant toric morphism $g \colon X\to X'$ to a divisorial
toric variety $X'$ and a rational map $h$ from $X'$ to $Y$ that is
defined on the image of $g$. By the universal properties of $p$ and
$q$, we have $g = g' \circ q \circ p$ with a toric morphism $g'
\colon Z \to X'$. So we arrive at the following commutative
diagram:
$$
\xymatrix{%
 X \ar[r]^{r} \ar[d]_{q \circ p}  & Y \ar[dl]|\alpha \\ 
Z \ar[r]_{g'} & X' \ar@{-->}[u]_{h} }
$$

Note that $g'(q(p(X))) = g(X)$ is contained in the domain
of definition of the rational map $h$. Since $r$ is surjective, we
have $q(p(X)) = \alpha(Y)$ and we obtain that $h$ is defined
on $g'(\alpha(Y))$. It follows that $(h \circ g') \circ \alpha$ 
is the identity on $Y$. This shows that $\alpha$ is injective.
Moreover, on the big torus of $Z$, the map $\alpha \circ (h \circ g')$
is the identity.

Consequently $\alpha \colon Y \to Z$ is a birational injection.
Since $Z$ is normal, Zariski's main theorem tells us that $\alpha$ is
in fact an open embedding. Since the image $\alpha(Y)$ is
invariant under the induced set-theoretical action of $T$ on $Y$,
the map $\alpha$ is an isomorphism. In particular, $r \colon X \to
Y$ is surjective. \endproof

We conclude this section with some examples. The above results 
in many situations give positive answers to the problem of existence
of quotients. A typical case are toric varieties defined by fans with
convex support:   

\begin{corollary}\label{convsupp1}
Let $X$ be  a toric variety arising from a fan with convex support. 
Then $X$ admits a divisorial reduction.
\end{corollary}

\proof Let the toric divisorial reduction $q \colon X \to X'$ 
arise from a map $Q \colon N \to N'$ of fans $\Delta$ and $\Delta'$.
Then $\sigma := Q(\vert \Delta \vert)$ is a convex cone in $N'$
and $\sigma\subset \vert \Delta' \vert$. 
Intersecting the cones of $\Delta'$ with $\sigma$, we obtain a 
further fan in $N'$, namely
$$ \Delta'' := \bigcup_{\tau' \in \Delta'}
                   \mathfrak{F} (\tau' \cap \sigma). $$

Let $X''$ be the associated toric variety. The identity map $N \to N'$
defines an affine toric morphism $g \colon X '' \to X'$. In particular, 
$X''$ is divisorial. Moreover, $Q \colon N \to N'$ is also a map of the
fans $\Delta$ and $\Delta''$. The corresponding toric
morphism $q' \colon X \to X''$ is surjective because $Q(\vert \Delta
\vert)$ equals $\vert \Delta''\vert$. Consider the decomposition
$$\xymatrix{
X \ar[rr]^{q} \ar[dr]_{q'} & & X' \\
& X'' \ar[ur]_{g}  } \,.$$
The universal property
of the toric divisorial reduction implies that $g \colon X'' \to X'$
is an isomorphism. Hence $q \colon X \to X'$ is surjective and the 
assertion follows from Corollary~\ref{divredchar}. \endproof

\begin{corollary}
Let $X$ be a divisorial toric variety arising from a fan with convex 
support. Then every subtorus action on $X$ admits a categorical 
quotient in the category of divisorial varieties.
\end{corollary}

\proof Let the toric quotient $p \colon X \to X'$ arise from a map
$P \colon N \to N'$ of fans $\Delta$ and $\Delta'$. 
By~\cite[Remark~2.5]{acha1}, each cone $\sigma' \in \Delta'$ is 
generated by images $P(\sigma)$ of certain $\sigma \in \Delta$.
Thus also $\Delta'$ has convex support and $p \colon X \to X'$
is surjective. So, Corollaries~\ref{divcatquot} and~\ref{convsupp1} give
the claim. \endproof

However, Corollary~\ref{divcatquot} also provides counterexamples to
existence of quotients. There can be different reasons for
nonsurjectivity of $q \circ p$, as the following examples show:

\begin{example}
For the toric variety $X$ described in Example~\ref{nichtsurjtdr}
the toric divisorial reduction is not surjective. Hence $X$ does not
admit a divisorial reduction. Moreover by Cox's construction,
see~\cite{Co}, $X$ is a good quotient of an open subset $\rq{X}
\subset \KK^{9}$ by a five dimensional subtorus $H \subset
(\KK^{*})^{9}$. So, the action of $H$ on $\rq{X}$ admits no
categorical quotient in the category of divisorial varieties. 
\end{example}

\begin{example}
Let $\Delta$ be the fan in $\ZZ^{4}$ having the following maximal
cones:
$$
\begin{array}{lcl}
\sigma_{1} := \cone((1,0,0,0),(0,1,0,0)), & & 
\sigma_{2} := \cone((0,0,1,0),(0,0,0,1))
\end{array}$$
The associated toric variety $X$ is an open toric subset of $\KK^{4}$.
Define a projection $P \colon \ZZ^{4} \to \ZZ^{3}$ by
$$\begin{array}{lcl}
P((1,0,0,0)) := (1,0,0), & & 
P((0,1,0,0)) := (0,1,0), \\ 
P((0,0,1,0)) := (0,0,1), & & 
P((0,0,0,1)) := (1,1,0).
\end{array}$$

By~\cite{acha1}, the toric morphism $p \colon X \to \KK^{3}$ 
defined by $P$ is the toric quotient for the action of
the subtorus $H := \ker(p)$ on $X$. Since $p$ is not surjective, the
action of $H$ on $X$ has no categorical quotient in the
category of divisorial varieties. 
\end{example}

\section{An open problem}\label{section7}

In this article we have solved the problem of existence of categorical
quotients for subtorus actions on toric varieties in the divisorial
category. For the analogous question in the category of all algebraic
varieties we have partial results.

For example, the toric quotient $p
\colon X \to X \tq H$ is a categorical quotient in the category of
algebraic varieties if the subtorus $H$ is of codimension at most two
\cite{acha4}, or if the map $p$ satisfies a certain curve lifting
property and $X \tq H$ is of expected dimension \cite{ac}. 

However, the general question still remains open. Therefore we pose it
here as a problem:

\begin{problem}
Give necessary and sufficient conditions for subtorus actions
on toric varieties to admit a categorical quotient in the category
of algebraic varieties.
\end{problem}

\bibliography{}

\end{document}